\documentclass[11pt]{amsart}
\usepackage{dsfont}
\usepackage{times}
\usepackage{enumerate}
\usepackage[euler-digits]{eulervm}
\theoremstyle{plain}
\usepackage{graphicx,color} 
\usepackage{mathtools}

\usepackage[hidelinks]{hyperref}
\usepackage{subcaption}
\usepackage{amsmath, amssymb, graphics}
\usepackage{float}

\newtheorem{theorem}{Theorem}[section] 
\newtheorem*{theorem*}{Theorem}
\newtheorem{cor}{Corollary}[section]

\newtheorem{prop}{Proposition}[section]

\newtheorem{remark}{Remark}[section]

\newtheorem*{conjecture*}{Conjecture}

\begin{document}
	
\title{Stability of Membranes}

\author{Bennett Palmer and \'Alvaro P\'ampano}

\date{\today}

\maketitle

\begin{abstract} In \cite{PP2}, the authors studied a particular class of equilibrium solutions of the Helfrich energy which satisfy a second order condition called the reduced membrane equation. In this paper we develop and apply a second variation formula for the Helfrich energy for this class of surfaces. 
	
The reduced membrane equation also arises as the Euler-Lagrange equation for the area of surfaces under the action of gravity in the three dimensional hyperbolic space. We study the second variation of this functional for a particular example.\\
	
\noindent K{\tiny EY} W{\tiny ORDS}.\, Helfrich Energy, Flux Formula, Second Variation.
	
\noindent MSC C{\tiny LASSIFICATION} (2020).\, 49Q10, 53A05.
\end{abstract}

\section{Introduction}

Introduced in 1973 by the German physicist Wolfgang Helfrich (\cite{H}), the energy functional
\begin{equation*}
	\mathcal{H}_{a,c_o,b}[\Sigma]:=\int_\Sigma \left(a \left[H+c_o\right]^2+b K\right)d\Sigma
\end{equation*}
is the principal tool for the determination of the morphology of biological membranes. In the  expression above, $H$ and $K$ are, respectively, the mean and Gaussian curvatures of the surface $\Sigma$ which models the membrane. The functional $\mathcal{H}_{a,c_o,b}$ contains three constant parameters, which depend on the composition of the membrane itself: the parameter $c_o$ is known as the spontaneous curvature (our definition of $c_o$ differs from the classical spontaneous curvature by the sign and a coefficient two) and the constants $a$ and $b$ are, respectively, the bending rigidity and saddle-splay moduli.

The Euler-Lagrange equation for the functional $\mathcal{H}_{a,c_o,b}$ is the fourth order nonlinear elliptic partial differential equation
\begin{equation*}
	\mathcal{L}[H+c_o]:=\Delta\left(H+c_o\right)+2\left(H+c_o\right)\left(H\left[H-c_o\right]-K\right)=0\,.
\end{equation*}

In a recent work (\cite{PP2}), the authors introduced a second order equation whose solutions automatically satisfy the above Euler-Lagrange equation $\mathcal{L}[H+c_o]=0$. The second order condition is referred to as the \emph{reduced membrane equation}. After a suitable normalization of the coordinates, this new equation takes the form
\begin{equation*}\label{RME1}
	H+c_o=\frac{\nu\cdot E}{A-X\cdot E}\,,
\end{equation*}
where $X$ is the position vector of the surface, $\nu$ is the unit surface normal, $E$ is a nonzero constant vector, and $A$ is a constant scalar. 

Conversely, the reduced membrane equation is a necessary condition for the Euler-Lagrange equation $\mathcal{L}[H+c_o]=0$ to hold if the surface is assumed to be an axially symmetric topological disc or topological sphere which is sufficiently regular so as to be a weak solution. However, there exist non-axially symmetric surfaces satisfying the reduced membrane equation as was shown in \cite{PP3}. Also, there is the well known family of examples, namely, the circular biconcave discoids (\cite{NOOY}) which, although they are axially symmetric topological spheres, do not satisfy the reduced membrane equation due to insufficient regularity at their poles. 

Via Noether's Theorem, each isometry of the three space generates a closed one form on a surface satisfying $\mathcal{H}_{a,c_o,b}$ and these one forms can be used to produce a  space of canonical solutions of  the equation $\mathcal{L}[f]=0$ corresponding to translations. The reduced membrane equation is, essentially, the condition that one of these solutions is a multiple of 
the function $H+c_o$.

If we make the further normalizations $E=E_3$, where $E_3$ is the constant unit vector field in the direction of the vertical coordinate $z$, and $A=0$, then the reduced membrane equation is itself the Euler-Lagrange equation for the functional
\begin{equation*}
	\mathcal{G}[\Sigma]:=\int_\Sigma \frac{1}{z^2}\,d\Sigma-2c_o\int_V \frac{1}{z^2}\,dV\,.
\end{equation*}
Here, we are assuming that $\Sigma$ is contained in one of the two open half-spaces determined by the plane $\{z=0\}$ and $V$ is the volume enclosed by the surface. When the background metric is changed to the hyperbolic one, then $\mathcal{G}$ represents the area of the surface in the hyperbolic space plus a constant times the gravitational potential energy corresponding to a vertical potential field. Therefore, the solutions of the reduced membrane equation can be interpreted as surfaces in the hyperbolic space whose mean curvature is a linear function of their height, i.e., they are capillary surfaces with gravity.

In this paper we develop a second variation formula for the Helfrich functional $\mathcal{H}_{a,c_o,b}$ in the limited case that the surface satisfies the reduced membrane equation. There are more general second variation formulas in the literature (see, for instance, \cite{CG,GTT}), but they are difficult to apply and their reduction to a simpler form when the reduced membrane equation holds is not at all clear. To circumvent this difficulty, we employ here a novel approach consisting on differentiating the flux formula. When the surface satisfies the reduced membrane equation, this method gives rise to a reasonable expression of the second variation formula. In this case, we show that the fourth order Jacobi operator for the Helfrich functional factors into a product of two second order operators, which involve the Jacobi operator for the functional $\mathcal{G}$. This is  analogous to the factorization of the Jacobi operator for the Willmore functional in the special case that the Willmore surface is conformal to a minimal surface in the three sphere (\cite{W}). 

Our second variation formula is applied to study the stability of solutions of the Euler-Helfrich functional. The Euler-Helfrich functional couples the Helfrich functional $\mathcal{H}_{a,c_o,b}$ on the interior of the surface to the elastic energy of the boundary and, hence, is given by
\begin{equation*}
\label{EH1}
E[\Sigma]:=\int_\Sigma \left(a[H+c_o]^2+bK\right)d\Sigma+\oint_{\partial \Sigma }\left(\alpha \kappa^2+\beta\right)ds\,,
\end{equation*}
where $\alpha,\beta>0$ are positive constants representing the flexural rigidity and the line tension, respectively. In the expression of $E$, $\kappa$ denotes the curvature of the boundary. This functional was introduced in a more general form in \cite{BMF} as a possible way to model bilipid membranes with edges. The stability analysis in this context is rather interesting since only stable equilibria are physically realizable.

We show that when $b=0$, a disc type surface satisfying the reduced membrane equation which is contained in a half-space $z\ge z_o>0$ and has its boundary in the plane $z=z_o$ is unstable. Also, for general values of $b$, a simple geometric condition is given which implies instability of an equilibrium surface.

In addition to studying the stability for the Euler-Helfrich functional, we also develop the variational theory of the functional $\mathcal{G}$. We use an approach developed by Maddocks (\cite{M}), Vogel (\cite{V}), and Koiso (\cite{K}), to use spectral methods for constrained variational problems. These results are applied to study the stability of a certain axially symmetric example from which a non symmetric branch of solutions bifurcates (\cite{PP3}). 

\section{Reduced Membrane Equation}

Throughout this paper $\Sigma$ will denote a compact (with or without boundary), connected and oriented surface which is immersed in the three dimensional Euclidean space through the $\mathcal{C}^4$ immersion 
$$X:\Sigma\longrightarrow{\bf R}^3\,.$$
We represent the canonical coordinates of ${\bf R}^3$ by $(x,y,z)$, while $E_i$, $i=1,2,3$, will denote the constant unit vector fields in the direction of the coordinate axes. The unit normal vector field $\nu:\Sigma\longrightarrow{\bf S}^2$ is defined by the property that it points out of any convex domain.

The \emph{Helfrich energy} for the immersion $X:\Sigma\longrightarrow{\bf R}^3$ is the functional
\begin{equation}\label{H}
	\mathcal{H}[\Sigma]:=\int_\Sigma \left(H+c_o\right)^2\,d\Sigma\,,
\end{equation}
where $c_o\in{\bf R}$ is any real constant and $H$ denotes the mean curvature of the immersion. When necessary, we will exhibit the dependence of the functional \eqref{H} on its parameter by writing $\mathcal{H}=\mathcal{H}_{c_o}$.

We will next review the first variation formula for the Helfrich energy $\mathcal{H}$. For details, the reader is referred to Section 2 and Appendix A of \cite{PP1}. Consider $\mathcal{C}^4$ variations of the immersion $X:\Sigma\longrightarrow{\bf R}^3$, i.e., $X+\epsilon\delta X+\mathcal{O}(\epsilon^2)$ where $\delta X$ is arbitrary, just assumed to be of class $\mathcal{C}^4$. To the variation field $\delta X$ we associate a vector field tangent to the surface by
\begin{equation}\label{J}
	\mathcal{J}_{\delta X}:=\left(H+c_o\right)\nabla\left[\nu\cdot\delta X\right]-\left(\nu\cdot\delta X\right)\nabla\left[H+c_o\right]+\left(H+c_o\right)^2\delta X^T\,.
\end{equation}
Then, the first variation of the Helfrich energy $\mathcal{H}$ can be expressed as
\begin{equation}\label{fvf}
	\delta_{\delta X}\mathcal{H}[\Sigma]=\int_\Sigma \mathcal{L}[H+c_o]\,\nu\cdot\delta X\,d\Sigma+\oint_{\partial\Sigma}\mathcal{J}_{\delta X}\cdot n\,ds\,,
\end{equation}
where $n$ is the outward pointing conormal to the boundary $\partial\Sigma$ (when the boundary is non-empty\footnote{Of course, in the case that the boundary $\partial\Sigma$ is empty, the second term in \eqref{fvf} does not appear.}) and $\mathcal{L}$ is the second order nonlinear elliptic operator defined by
\begin{equation}\label{mathcalL}
	\mathcal{L}[f]:=\Delta f+2\left(H[H-c_o]-K\right)f\,.
\end{equation}
Here, $K$ denotes the Gaussian curvature of the immersion.

\begin{remark} The first variation formula \eqref{fvf} of the Helfrich energy $\mathcal{H}$ is valid for any subdomain $\Omega\subset\Sigma$. 
\end{remark}

Considering compactly supported variations, it follows from the first variation formula \eqref{fvf} and the Fundamental Lemma of the Calculus of Variations that the Euler-Lagrange equation for the Helfrich energy $\mathcal{H}$ is the fourth order nonlinear elliptic partial differential equation
\begin{equation}\label{EL}
	\mathcal{L}[H+c_o]=\Delta\left(H+c_o\right)+2(H+c_o)\left(H[H-c_o]-K\right)=0\,.
\end{equation}
By expressing the equation $\mathcal{L}[H+c_o]=0$ in a nonparametric form we get a fourth order elliptic partial differential equation with coefficients that depend analytically on the height function and its derivatives. The real analyticity of any $\mathcal{C}^4$ critical surface then follows from elliptic regularity (see Theorem 6.6.1 of \cite{Morrey}).

In a previous paper (\cite{PP2}), by taking variation fields $\delta X\equiv E_i=\nabla X_i+\nu_i\nu$, $i=1,2,3$, we showed that any axially symmetric disc type surface which is critical for the Helfrich energy $\mathcal{H}$ is necessarily critical for a lower order functional. For the sake of completeness, we state this result here.

\begin{theorem}[\cite{PP2}, Theorem 4.1]\label{rme} Let $X:\Sigma\longrightarrow{\bf R}^3$ be an axially symmetric immersion of a topological disc critical for $\mathcal{H}$. Then, either the immersion has constant mean curvature, or else (after a suitable rigid motion and translation of the vertical coordinate) the equation
\begin{equation}\label{RME}
	H+c_o=-\frac{\nu_3}{z}\,,
\end{equation}
holds on $\Sigma$.
\end{theorem}

\begin{remark} The previous result can be extended to any surface that contains an axially symmetric topological disc on it. In particular, if the surface is a sufficiently regular axially symmetric topological sphere, the conclusion is still valid.
\end{remark}

The equation \eqref{RME} is referred to as the \emph{reduced membrane equation}. It turns out that any surface, whether it is axially symmetric or not and regardless of the topological type, on which \eqref{RME} holds, is also a solution of the Euler-Lagrange equation \eqref{EL}.  

\begin{prop}[\cite{PP2}, Proposition 4.1]\label{converse} Let $X:\Sigma\longrightarrow{\bf R}^3$ be an immersion whose mean curvature satisfies the reduced membrane equation \eqref{RME} on $\Sigma$. Then, the Euler-Lagrange equation $\mathcal{L}[H+c_o]=0$ given in \eqref{EL} also holds on $\Sigma$.
\end{prop}

\begin{remark}
	In the case where the spontaneous curvature $c_o$ is zero, the above result is more or less an analogue of the fact that minimal surfaces in space forms are Willmore surfaces. However, for $c_o\ne 0$ it is unexpected to see a similar result since the functional $\mathcal{H}$ is no longer conformally invariant.
\end{remark}

Additionally, the reduced membrane equation \eqref{RME} is the Euler-Lagrange equation for the functional (c.f., Theorem 4.2 and Remark 4.1 of \cite{PP2})
\begin{equation}\label{G}
	\mathcal{G}[\Sigma]:=\widetilde{\mathcal{A}}[\Sigma]-2c_o\int_{\widetilde{V}} \lvert z\rvert \,d\widetilde{V}=\widetilde{\mathcal{A}}[\Sigma]-2c_o\mathcal{U}[\Sigma]\,,
\end{equation}
where $\widetilde{\mathcal{A}}$ denotes the area of $\Sigma$, regarded as a surface in the hyperbolic space ${\bf H}^3$, $\widetilde{V}$ denotes the hyperbolic volume enclosed by $\Sigma$ and $\mathcal{U}$ is the gravitational potential energy of $\Sigma$, considered as a surface in ${\bf H}^3$. The functional $\mathcal{G}$ can be rewritten in terms of the metric on ${\bf R}^3$ as
$$\mathcal{G}[\Sigma]=\int_\Sigma \frac{1}{z^2}\,d\Sigma-2c_o\int_V \frac{1}{z^2}\,dV\,,$$
where $V$ is the Euclidean volume enclosed by $\Sigma$.

For convenience, we define a function $\xi:\Sigma\longrightarrow{\bf R}$ by
\begin{equation}\label{xi}
	\xi:=H+\frac{\nu_3}{z}\,.
\end{equation}
In terms of this new function $\xi$, the reduced membrane equation \eqref{RME} becomes $\xi=-c_o$.

Let $X:\Sigma\longrightarrow{\bf R}^3$ be an immersion satisfying the reduced membrane equation \eqref{RME} and consider a normal variation of this immersion given by $\delta X=f\nu$ where $f\in\mathcal{C}_o^\infty(\Sigma)$ is a compactly supported smooth function defined on $\Sigma$. We then define an operator $P$ using the pointwise variation of $\xi$ (as in \cite{PP3}). More precisely,
\begin{equation}\label{P}
	P[f]:=2\delta_{f\nu}\xi\,.
\end{equation}
The linearization of the reduced membrane equation \eqref{RME} for the critical immersion $X:\Sigma\longrightarrow{\bf R}^3$ is then $P[f]=0$.

\begin{prop}[\cite{PP3}] Let $X:\Sigma\longrightarrow{\bf R}^3$ be an immersion satisfying \eqref{RME} and assume that $X(\Sigma)$ is contained in one of the two open half-spaces determined by the plane $\{z=0\}$. Then, the operator $P$ defined in \eqref{P} can be expressed as
\begin{equation}\label{PL}
	P[f]=L[f]-2\frac{\nabla z\cdot\nabla f}{z}-2(H+c_o)^2f\,,
\end{equation}
where $f\in\mathcal{C}_o^\infty(\Sigma)$ and $L$ is the operator defined by the pointwise variation of the mean curvature for normal variations, namely, 
\begin{equation}\label{L}
	L[f]:=2\delta_{f\nu}H=\Delta f+\lVert d\nu\rVert^2f\,.
\end{equation}
(Recall that $\lVert d\nu\rVert^2=4H^2-2K$.)
\end{prop}

\begin{remark} Combining \eqref{mathcalL} and \eqref{PL}, the operator $P$ can as well be expressed as
\begin{equation}\label{Pf}
	P[f]=\mathcal{L}[f]-2\frac{\nabla z\cdot\nabla f}{z}-2c_o(H+c_o)f\,.
\end{equation}
In addition, a simple calculation using \eqref{mathcalL} and \eqref{Pf} shows that
\begin{equation}\label{easyP}
	\left(P+\frac{2}{z^2}\right)[f]=z \,\mathcal{L}\left[f/z\right].
\end{equation}
\end{remark}

Observe that the operator $P$ is not self-adjoint. Hence, for later use, we next compute the adjoint $P^*$ of the operator $P$.

\begin{prop} Let $X:\Sigma\longrightarrow{\bf R}^3$ be an immersion satisfying \eqref{RME} and assume that $X(\Sigma)$ is contained in one of the two open half-spaces determined by the plane $\{z=0\}$. Then, the adjoint operator of $P$ defined in \eqref{P} is given by
\begin{equation}\label{P*}
	P^*[f]:=L[f]+2\frac{\nabla z\cdot\nabla f}{z}-2\left(2H[H+c_o]+\frac{1}{z^2}\right)f\,,
\end{equation}
where $f\in\mathcal{C}_o^\infty(\Sigma)$ and $L$ is the operator given in \eqref{L}.
\end{prop}
{\it Proof.\:} Let $f,g\in\mathcal{C}_o^\infty(\Sigma)$ be two compactly supported smooth functions and consider the operator $P$ defined in \eqref{P}. From the definition of adjoint operator, we have that $P^*$ is characterized by
$$\int_\Sigma f P[g]\,d\Sigma=\int_\Sigma g P^*[f]\,d\Sigma\,.$$
Using the definition of $P$ given in \eqref{P} and \eqref{L}, we compute
$$\int_\Sigma f P[g]\,d\Sigma=\int_\Sigma f\left(\Delta g-2\frac{\nabla z\cdot\nabla g}{z}+\left[\lVert d\nu\rVert^2-2\left(H+c_o\right)^2\right]g\right)d\Sigma\,.$$
For the first term above we just integrate by parts twice, while for the second one integrating by parts, we have
$$\int_\Sigma f\frac{\nabla z\cdot\nabla g}{z}\,d\Sigma=\int_\Sigma\left(-g\frac{\nabla z\cdot\nabla f}{z}-fg\frac{\Delta z}{z}+fg\frac{\lVert \nabla z\rVert^2}{z^2}\right)d\Sigma\,.$$
Combining everything and using that $\lVert \nabla z\rVert^2=1-\nu_3^2$ and \eqref{RME},
$$\int_\Sigma f P[g]\,d\Sigma=\int_\Sigma g \left(\Delta f+2\frac{\nabla z\cdot\nabla f}{z}+\left[\lVert d\nu\rVert^2+2\frac{\Delta z}{z}-2\frac{1}{z^2}\right]f\right)d\Sigma\,.$$
Consequently, the adjoint operator $P^*$ is given by
\begin{equation*}P^*[f]=\Delta f+2\frac{\nabla z\cdot\nabla f}{z}+\left(\lVert d\nu\rVert^2+2\frac{\Delta z}{z}-2\frac{1}{z^2}\right)f\,.
\end{equation*}
Finally, we notice that $\Delta z=\Delta X\cdot E_3=2H\nu\cdot E_3=2H\nu_3$ and employ once again \eqref{RME} and \eqref{L} to conclude with the expression $P^*$ given in the statement. {\bf q.e.d.}

\begin{remark}\label{rem} Combining \eqref{mathcalL} and \eqref{P*}, the operator $P^*$ can as well be expressed as
\begin{equation}\label{P*f}
	P^*[f]=\mathcal{L}[f]+2\frac{\nabla z\cdot\nabla f}{z}-2\left(H[H+c_o]+\frac{1}{z^2}\right)f\,.
\end{equation}
In addition, a simple calculation using \eqref{mathcalL} and \eqref{P*f} shows that
\begin{equation}\label{easyP*}
	\left(P^*+\frac{2}{z^2}\right)[f]=\frac{1}{z} \mathcal{L}[z f]\,.
\end{equation}
\end{remark}

We finish this section by proving a couple of relations between the operator $P$ and its adjoint $P^*$. The first one is an adjoint formula for smooth functions that may not vanish on the boundary.

\begin{prop} Let $X:\Sigma\longrightarrow{\bf R}^3$ be an immersion satisfying \eqref{RME} and assume that $X(\Sigma)$ is contained in one of the two open half-spaces determined by the plane $\{z=0\}$. Then, 
\begin{equation}\label{int}
	\int_\Sigma\left(g\, P[f]-f\, P^*[g]\right)d\Sigma=\oint_{\partial\Sigma}\left(g\, \partial_n f-f\,\partial_n g-2fg\,\frac{\partial_nz}{z}\right)ds\,,
\end{equation}
holds for every $f,g\in\mathcal{C}^\infty(\Sigma)$.
\end{prop}
{\it Proof.\:} Let $P$ and $P^*$ be the operators given in \eqref{PL} and \eqref{P*}, respectively. Using these expressions, we compute for every $f,g\in\mathcal{C}^\infty(\Sigma)$,
\begin{eqnarray*}
	\int_\Sigma \left(gP[f]-fP^*[g]\right)d\Sigma&=&\int_\Sigma\left(g\Delta f-f\Delta g\right)d\Sigma\\
	&+&\int_\Sigma\left(2\left[H^2-c_o^2+\frac{1}{z^2}\right]fg-2\frac{\nabla z\cdot\nabla(fg)}{z}\right)d\Sigma\,.
\end{eqnarray*}
Observe that since $\Delta z=2H\nu_3$ and \eqref{RME} holds,
$$\nabla\cdot\left(\frac{\nabla z}{z}\right)=\frac{\Delta z}{z}-\frac{\lVert \nabla z\rVert^2}{z^2}=\frac{\Delta z}{z}-\frac{1}{z^2}+\frac{\nu_3^2}{z^2}=-\left(H^2-c_o^2+\frac{1}{z^2}\right),$$
and, hence, the integrand on the second line above is
\begin{eqnarray*}
-2\nabla\cdot\left(\frac{\nabla z}{z}\,fg\right)&=&-2\nabla\cdot\left(\frac{\nabla z}{z}\right)fg-2\frac{\nabla z\cdot\nabla(fg)}{z}\\
&=&2\left(H^2-c_o^2+\frac{1}{z^2}\right)fg-2\frac{\nabla z\cdot\nabla(fg)}{z}\,.
\end{eqnarray*}
Consequently,
\begin{eqnarray*}
	\int_\Sigma \left(gP[f]-fP^*[g]\right)d\Sigma&=&\int_\Sigma\left(g\Delta f-f\Delta g-2\nabla\cdot\left[\frac{\nabla z}{z}\,fg\right]\right)d\Sigma\\
	&=&\oint_{\partial\Sigma}\left(g\,\partial_n f-f\,\partial_ng-2fg\frac{\partial_nz}{z}\right)ds\,,
\end{eqnarray*}
where for the first two terms we have applied Green's second identity, while for the last one we have used the Divergence Theorem. {\bf q.e.d.}
\\

The second relation shows that for suitable smooth functions the operator $P^*$ at those functions can be computed in terms of the operator $P$.

\begin{prop} Let $X:\Sigma\longrightarrow{\bf R}^3$ be an immersion satisfying \eqref{RME} and assume that $X(\Sigma)$ is contained in one of the two open half-spaces determined by the plane $\{z=0\}$. Then, 
\begin{equation}\label{PP*}
	P^*\left[\frac{f}{z^2}\right]=\frac{1}{z^2}P[f]\,,
\end{equation}
holds for every $f\in\mathcal{C}^\infty(\Sigma)$.
\end{prop}
{\it Proof.\:} Let $g\in\mathcal{C}^\infty(\Sigma)$. We evaluate \eqref{easyP*} at $f=g/z^2$ to obtain
$$P^*\left[\frac{g}{z^2}\right]+2\,\frac{g}{z^4}=\frac{1}{z}\mathcal{L}\left[g/z\right].$$
We now use \eqref{easyP} to substitute $\mathcal{L}[g/z]$ in above expression, concluding that
$$P^*\left[\frac{g}{z^2}\right]+2\,\frac{g}{z^4}=\frac{1}{z^2}\left(P[g]+2\,\frac{g}{z^2}\right).$$
The second terms on the left and right hand sides cancel out and, then, we end up with the expression of the statement. {\bf q.e.d.}

\section{Second Variation Through the Flux Formula}

In this section we will develop some of the variational theory of the Helfrich energy $\mathcal{H}$, concentrating on surfaces satisfying the reduced membrane equation \eqref{RME}. The operators $P$ and its adjoint $P^*$ defined in previous section will play an essential role in the second variation of the Helfrich energy $\mathcal{H}$.

A direct calculation of the second variation of $\mathcal{H}$ is problematic since the expression involved is very complicated, which makes it almost impossible to deal with. Hence, we carry out here the calculation using the flux formula. 

\subsection{Area Functional} In order to demonstrate the method, we begin with the far simpler calculation of the second variation of the area functional for minimal surfaces. For an immersion $X:\Sigma\longrightarrow{\bf R}^3$ consider the area functional
\begin{equation}\label{area}
	\mathcal{A}[\Sigma]:=\int_\Sigma d\Sigma\,,
\end{equation}
and arbitrary variations of the immersion, that is, deformations $X+\epsilon\delta X+\mathcal{O}(\epsilon^2)$. Then, the first variation of $\mathcal{A}$ is given by
\begin{equation}\label{fvfA}
	\delta_{\delta X}\mathcal{A}[\Sigma]=-2\int_\Sigma H\,\nu\cdot \delta X\,d\Sigma+\oint_{\partial\Sigma} n\cdot \delta X\,ds\,.
\end{equation}
Critical points of the area functional are characterized by $H\equiv 0$ on the interior of $\Sigma$. In other words, they are minimal immersions.

Take a variation field $\delta X\equiv E_i$, $i=1,2,3$. Since $E_i$ is a constant vector field which generates a symmetry (more precisely, a translation) of the Lagrangian, it must be the case that $\delta_{E_i}\mathcal{A}[\Omega]=0$ for any arbitrary smoothly bounded relatively compact subdomain $\Omega\subset\Sigma$ and, hence, from \eqref{fvfA} we have
\begin{equation}\label{cl}
	2\int_\Omega H\nu_i\,d\Sigma=\oint_{\partial\Omega} n_i\,ds\,,
\end{equation}
for every $i=1,2,3$.

We will next linearize the previous equation at a minimal immersion $X:\Sigma\longrightarrow{\bf R}^3$ in order to obtain the variation of the mean curvature. Consider an arbitrary normal variation given by the variation field $\delta X=f\nu$, $f\in\mathcal{C}_o^\infty(\Sigma)$. Differentiating \eqref{cl}, we obtain
\begin{equation}\label{delta}
	2\int_\Omega\left(\delta_{f\nu} H\right)\nu_i\,d\Sigma=\delta_{f\nu}\left(\oint_{\partial\Omega}n_i\,ds\right)=\delta_{f\nu}\left(\oint_{\partial\Omega}\star dx_i\right),
\end{equation}
since $H=0$ holds on $\Omega$. In the second equality above we are using an alternative expression involving the codifferential of a function, which will be more convenient for differentiation purposes (c.f., Remark \ref{alternative}). Indeed, with respect to the normal variation $\delta X=f\nu$, the variation of the codifferential of $x_i$ is given by
$$\delta\left(\star d\right)x_i=-2f\star \left(d\nu(\nabla x_i)\right)^\flat=-2f\star d\nu_i\,,$$
where $\flat$ denotes the flat of a vector field. Substituting this formula in \eqref{delta}, we compute
\begin{eqnarray*}
	2\int_\Omega \left(\delta_{f\nu} H\right)\nu_id\Sigma&=&\oint_{\partial\Omega}\left(\star d(\delta_{f\nu} x_i)+\delta_{f\nu}(\star d)x_i\right)\\
	&=&\oint_{\partial\Omega}\left(\star d(f\nu_i)-2f\star d\nu_i\right)=\oint_{\partial\Omega}\left(\nu_i\star df-f\star d\nu_i\right)\\
	&=&\int_\Omega \left(\nu_i\Delta f-f \Delta \nu_i\right)d\Sigma\,,
\end{eqnarray*}
where in the last equality we have used Green's second identity. Moreover, for a minimal surface $\Delta\nu_i=2K\nu_i$ holds and, hence,
$$2\int_\Omega \left(\delta_{f\nu} H\right)\nu_i\,d\Sigma=\int_\Omega\left(\Delta f-2Kf\right)\nu_i\,d\Sigma\,.$$
Therefore, since $\Omega$ is an arbitrary subdomain of $\Sigma$, we conclude that the Jacobi operator for the area functional $\mathcal{A}$ is
$$2\delta_{f\nu} H=\Delta f-2Kf\,,$$
as expected (this is, precisely, the operator $L$ given in \eqref{L} when $H\equiv 0$ holds, i.e., for minimal immersions).

\begin{remark}\label{alternative} The same computation can be carried out without employing the codifferential $\star d x_i$. Indeed, using that $\delta n=\tau_g f T+\partial_nf \nu$ for $\delta X=f\nu$ (cf. Appendix of \cite{PP0}), where $\{n,T,\nu\}$ is the Darboux frame along $\partial\Omega$ and $\tau_g$ is the geodesic torsion of $\partial\Omega$, we get from \eqref{delta} that
\begin{eqnarray*}
	2\int_\Omega \left(\delta_{f\nu} H\right)\nu_i\,d\Sigma&=&\delta_{f\nu}\left(\oint_{\partial\Omega}n_i\,ds\right)=\oint_{\partial\Omega}\left(\delta_{f\nu} n_i\,ds+n_i \delta_{f\nu}(ds)\right)\\
	&=&\oint_{\partial\Omega}\left(\tau_g f\, T_i+\partial_n f \,\nu_i-\kappa_n f\, n_i\right)ds\\&=&\oint_{\partial\Omega}\left(\partial_n f\, \nu_i-f\, \partial_n \nu_i\right)ds\,,
\end{eqnarray*}
since $\partial_n\nu_i=d\nu(n)_i=\kappa_n n_i-\tau_g T_i$, where $\kappa_n$ is the normal curvature of $\partial\Omega$. Then, applying Green's second identity and proceeding in the same way as above we conclude with the result.
\end{remark}

We finish the demonstration of the method by employing the Jacobi operator $L$ computed above to obtain the second variation of the area functional $\mathcal{A}$ for minimal immersions. The following result is classical so we simply state it here without a proof.

\begin{prop} Let $X:\Sigma\longrightarrow{\bf R}^3$ be a minimal immersion. Then, 
	$$\delta^2_{f\nu}\mathcal{A}[\Sigma]=-\int_\Sigma f\left(\Delta f-2Kf\right)d\Sigma=-\int_\Sigma f\,L[f]\,d\Sigma\,,$$
holds for every $f\in\mathcal{C}_o^\infty(\Sigma)$.
\end{prop}

\subsection{Helfrich Energy} We now apply the idea of computing the second variation through the flux formula to the more complicated case of the Helfrich energy $\mathcal{H}$.

Let $\Omega$ be a smoothly bounded relatively compact subdomain of $\Sigma$ and take the variation vector field $\delta X=E_3$, which generates a symmetry of the Lagrangian and, hence, $\delta_{E_3}\mathcal{H}[\Omega]=0$ must hold\footnote{Similar computations work for the vector fields $E_i$, $i=1,2$. However, up to a rigid motion, we may simply consider $E_3$.}. Then, from the first variation of the Helfrich energy $\mathcal{H}$ given in \eqref{fvf}, we have
\begin{equation}\label{deltaH}
	\int_\Omega \mathcal{L}[H+c_o]\nu_3\,d\Sigma=-\oint_{\partial\Omega} (H+c_o)^2\partial_n\left(\frac{\nu_3}{H+c_o}+z\right)ds\,.
\end{equation}

Assume now that $X:\Sigma\longrightarrow{\bf R}^3$ is an immersion satisfying the reduced membrane equation \eqref{RME}. From Proposition \ref{converse} it then follows that \eqref{EL} also holds on $\Sigma$, that is, $\mathcal{L}[H+c_o]=0$. Moreover, observe that if \eqref{RME} holds, then $$\partial_nH=-\frac{\partial_n\nu_3}{z}+\nu_3\frac{\partial_nz}{z^2}\,,$$ 
and
$$\partial_n\left(\frac{\nu_3}{H+c_o}+z\right)=\frac{\partial_n\nu_3}{H+c_o}-\frac{\nu_3\partial_nH}{(H+c_o)^2}+\partial_nz=0\,.$$
Therefore, for a normal variation given by the variation field $\delta X=f\nu$, $f\in\mathcal{C}_o^\infty(\Sigma)$, differentiating equation \eqref{deltaH}, we get
\begin{equation}\label{deltaH2}
	\int_\Omega \delta_{f\nu}(\mathcal{L}[H+c_o])\nu_3d\Sigma=-\oint_{\partial\Omega}(H+c_o)^2\partial_n\left(\delta_{f\nu}\left[\frac{\nu_3}{H+c_o}+z\right]\right)ds.
\end{equation}
Taking into account that \eqref{RME} is equivalent to $\xi=-c_o$ where $\xi$ is the function defined in \eqref{xi}, we next compute
\begin{equation*}
	\delta_{f\nu}\left(\frac{\nu_3}{H+c_o}+z\right)=\delta_{f\nu}\left(\frac{[\xi+c_o]z}{H+c_o}\right)=\delta_{f\nu}\xi\,\frac{z}{H+c_o}=\frac{1}{2}\frac{z P[f]}{H+c_o}\,,
\end{equation*}
where we have used the definition of the operator $P$ given in \eqref{P}. Now, expanding the normal derivative, we get
$$\partial_n\left(\delta_{f\nu}\left[\frac{\nu_3}{H+c_o}+z\right]\right)=\frac{1}{2}\left(\frac{\partial_n\left(z P[f]\right)}{H+c_o}-\frac{z P[f] \partial_n(H+c_o)}{(H+c_o)^2}\right).$$
Therefore, multiplying this by $(H+c_o)^2$, we obtain that the right hand side of \eqref{deltaH2} is
$$-\frac{1}{2}\oint_{\partial\Omega} \left([H+c_o]\partial_n\left( z P[f]\right)-z P[f]\partial_n(H+c_o)\right)ds\,.$$
Applying Green's second identity, this is the same as
$$-\frac{1}{2}\int_\Omega \left([H+c_o]\Delta\left(z P[f]\right)-z P[f]\Delta(H+c_o)\right)d\Sigma\,.$$
We now use that $\mathcal{L}[H+c_o]=0$ to substitute $\Delta(H+c_o)$ in the above expression and the definition of the operator $\mathcal{L}$ given in \eqref{mathcalL}, to obtain
\begin{eqnarray*}
	\int_\Omega\delta_{f\nu}\left(\mathcal{L}[H+c_o]\right)\nu_3\,d\Sigma&=&-\frac{1}{2}\int_\Omega (H+c_o)\mathcal{L}\left[z P[f]\right]\,d\Sigma\\
	&=&\frac{1}{2}\int_\Omega \frac{1}{z}\mathcal{L}\left[z P[f]\right]\nu_3\,d\Sigma\,,
\end{eqnarray*}
where we have used \eqref{RME} again.

Consequently, since $\Omega$ is an arbitrary subdomain of $\Sigma$, we have that 
$$\delta_{f\nu}\left(\mathcal{L}[H+c_o]\right)=\frac{1}{2z}\mathcal{L}\left[zP[f]\right]\,.$$ 
Finally, we use the relation \eqref{easyP*} of Remark \ref{rem} between the operators $\mathcal{L}$ and $P^*$ to deduce that the linearization of $\mathcal{L}[H+c_o]$ is the self-adjoint operator $F$ defined by
\begin{equation}\label{F}
	F[f]:=\delta_{f\nu}\left(\mathcal{L}[H+c_o]\right)=\frac{1}{2}\left(P^*+\frac{2}{z^2}\right)\circ P[f]\,,
\end{equation}
where $\circ$ denotes the composition of operators.

\begin{remark} The self-adjoint operator $F$ given in \eqref{F} is regular since it arises from the variation of a smooth quantity. However, if $z=0$ somewhere on the surface, the expressions of $P[f]$ and $P^*[f]$, \eqref{PL} and \eqref{P*} respectively, may have singularities even if $f$ is smooth. This makes sense from the definition of $P$ given in \eqref{P}. Indeed, the quantity $\xi=H+\nu_3/z$ may not be regular if the surface is perturbed.
\end{remark}

Using the operator $F$ we can now compute the second variation of the Helfrich energy $\mathcal{H}$ for surfaces satisfying the reduced membrane equation \eqref{RME}.

\begin{theorem}\label{second} Let $X:\Sigma\longrightarrow{\bf R}^3$ be an immersion critical for $\mathcal{H}$ satisfying the reduced membrane equation \eqref{RME}. Then, for every $f\in\mathcal{C}_o^\infty(\Sigma)$,
	\begin{equation}\label{svf1}
		\delta_{f\nu}^2\mathcal{H}[\Sigma]=\int_\Sigma f\,F[f]\,d\Sigma+\frac{1}{2}\oint_{\partial\Sigma}L[f]\,\partial_nf\,ds\,,
	\end{equation}
where $F$ is the fourth order self-adjoint operator defined in \eqref{F} and $L$ is the second order operator defined in \eqref{L}. Moreover, if $X(\Sigma)$ is contained in one of the two open half-spaces determined by the plane $\{z=0\}$, \eqref{svf1} can be rewritten as
	\begin{equation}\label{svf2}
		\delta_{f\nu}^2\mathcal{H}[\Sigma]=\frac{1}{2}\int_\Sigma P[f]\left(P+\frac{2}{z^2}\right)[f]\,d\Sigma+\oint_{\partial\Sigma}\left(\partial_nf\right)^2\frac{\partial_nz}{z}\,ds\,,
	\end{equation}
where $P$ is the second order operator defined in \eqref{P}.
\end{theorem}
{\it Proof.\:} Let $X:\Sigma\longrightarrow{\bf R}^3$ be an immersion satisfying the reduced membrane equation \eqref{RME} and consider arbitrary compactly supported normal variations of $X$ given by the variation vector field $\delta X=f\nu$, $f\in\mathcal{C}_o^\infty(\Sigma)$.

Then, the first variation of the Helfrich energy \eqref{fvf} reduces to
$$\delta_{f\nu}\mathcal{H}[\Sigma]=\int_\Sigma\mathcal{L}[H+c_o]f\,d\Sigma+\oint_{\partial\Sigma}(H+c_o)\,\partial_n f\,ds\,,$$
since the first term in $\mathcal{J}_{\delta X}$, \eqref{J}, vanishes because $f=0$ on $\partial\Sigma$, and the last term is zero because the variations are normal. From this it follows that if the immersion $X:\Sigma\longrightarrow{\bf R}^3$ is critical for $\mathcal{H}$ then, in addition to $\mathcal{L}[H+c_o]=0$ (which is satisfied since \eqref{RME} holds, c.f., Proposition \ref{converse}), the boundary condition $H+c_o=0$ must also hold.

Assume that $X:\Sigma\longrightarrow{\bf R}^3$ is a critical immersion satisfying \eqref{RME}. Then, differentiating the first variation of $\mathcal{H}$, we get
$$\delta_{f\nu}^2\mathcal{H}[\Sigma]=\int_\Sigma \delta_{f\nu}\left(\mathcal{L}[H+c_o]\right)f\,d\Sigma+\oint_{\partial\Sigma}\left(\delta_{f\nu}H\right)\partial_n f\,ds\,.$$
Finally, from the definitions of the operators $F$ and $L$, \eqref{F} and \eqref{L} respectively, we conclude with \eqref{svf1}. This finishes the proof of the first part of the statement.

For the second part, assume also that $X(\Sigma)$ is contained in one of the two half-spaces determined by the plane $\{z=0\}$. Then, for every $f,g\in\mathcal{C}^\infty(\Sigma)$, the operators $P$ and $P^*$ are related by \eqref{int}. Taking $f\in\mathcal{C}_o^\infty(\Sigma)$ (i.e., $f=0$ on $\partial\Sigma$) and $g=P[f]$, \eqref{int} gives
$$\int_\Sigma\left(\left(P[f]\right)^2-f\,P^*\left[P[f]\right]\right)d\Sigma=\oint_{\partial\Sigma}P[f]\,\partial_nf\,ds\,.$$
Substituting this in \eqref{svf1} and using the definition of $F$ given in \eqref{F}, we obtain
\begin{eqnarray*}
	\delta_{f\nu}^2\mathcal{H}[\Sigma]&=&\frac{1}{2}\int_\Sigma\left(P^*\left[P[f]\right]+\frac{2P[f]}{z^2}\right)f\,d\Sigma+\frac{1}{2}\oint_{\partial\Sigma}L[f]\,\partial_nf\,ds\\
	&=&\frac{1}{2}\int_\Sigma\left(\left(P[f]\right)^2+\frac{2fP[f]}{z^2}\right)d\Sigma+\frac{1}{2}\oint_{\partial\Sigma}\left(L[f]-P[f]\right)\partial_nf\,ds\\
	&=&\frac{1}{2}\int_\Sigma P[f]\left(P+\frac{2}{z^2}\right)[f]\,d\Sigma+\oint_{\partial\Sigma}\frac{\nabla z\cdot\nabla f}{z}\,\partial_nf\,ds\,,
\end{eqnarray*}
where in the last equality we have used the expression of $P$ given in \eqref{PL} and that $f=0$ holds on $\partial\Sigma$. Moreover, since $f=0$ on $\partial\Sigma$, $\partial_sf=0$ holds too and $\nabla f=\partial_n f n$, drawing the desired conclusion. {\bf q.e.d.}

\begin{remark} Observe that the fourth order operator $F$ arising in the second variation of the Helfrich energy \eqref{svf1} is defined as the composition of two second order operators, namely, $P^*+2/z^2$ and $P$. However, in \eqref{svf2} we have rewritten the second variation of the Helfrich energy as the product of two second order operators, $P$ and $P+2/z^2$. 

In the paper \cite{W}, Weiner studied the stability of closed minimal surfaces for the Willmore functional $\mathcal{H}_o$, by factoring its Jacobi operator as $L(L-2)$, where $L$ is the Jacobi operator for the area (see \eqref{L}). He concluded that the surface is stable if the spectrum of $L$ has no eigenvalues in the interval $(-2,0)$. Here, we have an analogous situation where the Jacobi operator $F$ factorizes involving the linearization of the reduced membrane equation \eqref{RME}, namely, $P$.
\end{remark}

In the following result we relate an eigenvalue problem for the operator $P$ with that of $F$. Specifically, we consider the problem
\begin{equation}\label{EigenType}
	P[f]+\frac{\lambda}{z^2}f=0\,,
\end{equation}
with $f$ vanishing on $\partial \Sigma$ if this boundary is non-empty.

\begin{prop}\label{FP} Let $X:\Sigma\longrightarrow{\bf R}^3$ be an immersion satisfying \eqref{RME} and assume that $X(\Sigma)$ is contained in one of the two open half-spaces determined by the plane $\{z=0\}$. If the function $f$ satisfies \eqref{EigenType} for some $\lambda\in{\bf R}$. Then,
	\begin{equation}\label{Flambda}
		F[f]=\frac{1}{2}\,\lambda\left(\lambda-2\right)\frac{f}{z^4}\,,
	\end{equation}
	holds.
\end{prop}
{\it Proof.\:} Assume that the second order operator $P$ satisfies $P[f]=-\lambda f/z^2$. It then follows from the definition of the operator $F$ given in \eqref{F} and the relation \eqref{PP*} between $P$ and $P^*$ that
\begin{eqnarray*}
	F[f]&=&\frac{1}{2}\left(P^*+\frac{2}{z^2}\right)\left[P[f]\right]=\frac{1}{2}\left(P^*+\frac{2}{z^2}\right)\left[\frac{-\lambda f}{z^2}\right]\\ &=&-\frac{\lambda}{2}\left(P^*\left[\frac{f}{z^2}\right]+2\frac{f}{z^4}\right)=-\frac{\lambda}{2z^2}\left(P[f]+2\frac{f}{z^2}\right)\\
	&=&\frac{1}{2}\,\lambda\left(\lambda-2\right)\frac{f}{z^4}\,,
\end{eqnarray*}
proving the result. {\bf q.e.d.}

\begin{remark}
	The operator $P$ is self-adjoint with respect to the measure $z^{-2}d\Sigma$, so the spectrum can be determined in the usual way from the Rayleigh quotient
	\begin{equation*}\label{RQ} 
		\frac{-\int_\Sigma f\,P[f]\:z^{-2}d\Sigma}{\int_\Sigma f^2z^{-4}d\Sigma}\:.
	\end{equation*}

	In \cite{PP3}, using variations $X_\rho:=X+\rho E_3$ and the definition of the operator $P$ given in \eqref{P}, it was shown that
	\begin{equation}\label{Pnu3}
		P[\nu_3]+2\frac{\nu_3}{z^2}=0\,.
	\end{equation}
	Therefore, when $f=\nu_3$, we obtain the value $\lambda=2$ in \eqref{EigenType}. Moreover, similar computations using variations $X+\rho E_i$, $i=1,2$, show that when $f=\nu_i$, $i=1,2$, we get the value $\lambda=0$ in \eqref{EigenType}. Consequently, we deduce from \eqref{Flambda} that $F[\nu_i]=0$ for every $i=1,2,3$.
\end{remark}

\section{Applications to the Euler-Helfrich Functional}

In this section we will apply our findings about the second variation of the Helfrich energy $\mathcal{H}$ to study the stability of equilibria for the Euler-Helfrich functional. For a smooth immersion $X:\Sigma\longrightarrow{\bf R}^3$of a compact surface $\Sigma$ with smooth boundary, the \emph{Euler-Helfrich energy} functional is the total energy
\begin{equation}\label{E}
E[\Sigma]:=\int_\Sigma \left(a[H+c_o]^2+bK\right)d\Sigma+\oint_{\partial \Sigma }\left(\alpha \kappa^2+\beta\right)ds\,,
\end{equation}
where $a>0$, $\alpha>0$, $\beta>0$ and $c_o,b$ are any real constants. Here, $\kappa$ denotes the (Frenet) curvature of the boundary $\partial\Sigma$. As customary, we will write $E=E_{a,c_o,b,\alpha,\beta}$ when we need to explicitly exhibit the dependence of \eqref{E} on its parameters.

In order to describe the equations characterizing the equilibria of $E$, we define the vector field along $\partial\Sigma$
\begin{equation}\label{Jboundary}
 J:=2\alpha T''+\left(3\alpha\kappa^2-\beta\right)T\,,
\end{equation}
where $T$ is the unit tangent field to $\partial\Sigma$ and $\left(\,\right)'$ denotes the derivative with respect to $s$, the arc length parameter of the boundary. Using this and \eqref{J}, we get that the first variation of $E$ is (see \cite{PP1,PP2} for details)
\begin{eqnarray}\label{var1}
\delta_{\delta X}E[\Sigma]&=& a\int_\Sigma \mathcal{L}[H+c_o]\nu\cdot\delta X\:d\Sigma+\oint_{\partial \Sigma} \left([a\mathcal{J}_{\delta X}+bK\delta X]\cdot n+J'\cdot \delta X\right)ds\nonumber\\
&&+\oint_{\partial\Sigma}b\left(\left[d\nu+2H\,{\rm Id}\right]\nabla\left[\nu\cdot\delta X\right]\right)\cdot n\,ds\:.   
\end{eqnarray}
Therefore, the Euler-Lagrange equations associated to $E$ are
\begin{eqnarray}
\mathcal{L}[H+c_o]&=&0\,,\quad\quad\quad\text{on $\Sigma$}\,,\label{EL1}\\
a\left(H+c_o\right)+b\kappa_n&=&0\,,\quad\quad\quad\text{on $\partial\Sigma$}\,,\label{EL2}\\
J'\cdot\nu-a\partial_n H+b\tau_g'&=&0\,,\quad\quad\quad\text{on $\partial\Sigma$}\,,\label{EL3}\\
J'\cdot n+a\left(H+c_o\right)^2+bK&=&0\,,\quad\quad\quad\text{on $\partial\Sigma$}\,,\label{EL4}
\end{eqnarray}
where $\mathcal{L}[H+c_o]$ is given in \eqref{EL}.

We next compute the second variation of $E$ for compactly supported normal variations. 

\begin{theorem}\label{2E}
	Let $X:\Sigma\longrightarrow{\bf R}^3$ be an immersion critical for $E$ satisfying the reduced membrane equation \eqref{RME} and assume that $X(\Sigma)$ is contained in one of the two open half-spaces determined by the plane $\{z=0\}$. Then, 
	\begin{equation}\label{2vfE}
		\delta_{f\nu}^2E[\Sigma]=\frac{a}{2}\int_\Sigma P[f]\left(P+\frac{2}{z^2}\right)[f]\,d\Sigma+\oint_{\partial\Sigma} \left(a\frac{\partial_n z}{z}-b\kappa_g\right)\left(\partial_nf\right)^2ds\,,
	\end{equation}
for every $f\in\mathcal{C}_o^\infty(\Sigma)$.
\end{theorem}
{\it Proof.\:} Let $X:\Sigma\longrightarrow{\bf R}^3$ be an immersion critical for $E$ satisfying \eqref{RME} and consider arbitrary compactly supported normal variations given by $\delta X=f\nu$, $f\in\mathcal{C}_o^\infty(\Sigma)$.

From the definitions of the vector fields $\mathcal{J}_{\delta X}$ and $J$, \eqref{J} and \eqref{Jboundary} respectively, we obtain that the first variation of $E$ given in \eqref{var1} reduces to
$$\delta_{f\nu}E[\Sigma]=a\int_\Sigma\mathcal{L}[H+c_o]f\,d\Sigma+\oint_{\partial\Sigma}\left(a[H+c_o]+b\kappa_n\right)\partial_nf\,ds\,,$$
since $f=0$ on $\partial\Sigma$.

Therefore, computing the variation of this quantity and taking into account that \eqref{EL1} and \eqref{EL2} hold, we obtain
$$\delta_{f\nu}^2E[\Sigma]=a\int_\Sigma \delta(\mathcal{L}[H+c_o])f\,d\Sigma+\oint_{\partial\Sigma}\left(a\delta_{f\nu}H+b\delta_{f\nu}[\kappa_n]\right)\partial_n f\,ds\,.$$
Observe that the first two terms are a multiple of $\delta_{f\nu}^2\mathcal{H}[\Sigma]$ (c.f., Theorem \ref{second}) and so they can be rewritten as a suitable multiple of the right hand side of \eqref{svf2}. It remains to compute the last term. 

Since both $f$ and $\partial_sf$ vanish on $\partial\Sigma$, we have $\partial_s(\delta X)=\partial_s(f\nu)=0$ and, hence, in our case, for any smooth function $g$ it holds along $\partial\Sigma$
$$\delta\left(\partial_sg\right)=\partial_s\left(\delta g\right)-\left(\partial_s[\delta X]\cdot T\right)\partial_sg\,.$$

Therefore, using this and from the definition of the normal curvature of $\partial\Sigma$, namely, $\kappa_n=-T\cdot\partial_s\nu$ (c.f., \cite{PP2}), we conclude that
\begin{equation*}
	\delta_{f\nu}\left(\kappa_n\right)=-\delta_{f\nu}(T)\cdot \partial_s\nu-T\cdot\delta_{f\nu}(\partial_s\nu)=-T\cdot\partial_s(\delta_{f\nu}\nu)=\nabla_T\left(\nabla f\right)\cdot T\,.
\end{equation*}
However, $\nabla f=\partial_nf n$ on $\partial\Sigma$ since $\partial_sf=0$ and, hence,
$$\delta_{f\nu}\left(\kappa_n\right)=\partial_n f\,\partial_s n\cdot T=-\kappa_g\,\partial_nf\,.$$
Finally, combining this computation and \eqref{svf2}, we prove the result. {\bf q.e.d.}
\\

As a consequence of the second variation of the Euler-Helfrich functional obtained in Theorem \ref{2E} we have the following result.

\begin{cor}\label{cor} Let $X:\Sigma\longrightarrow{\bf R}^3$ be an immersion of a disc type surface critical for $E$ with $c_o>0$ and $b=0$ satisfying the reduced membrane equation \eqref{RME}. If $X(\Sigma)$ is contained in $\{z\geq z_o>0\}$ and the boundary $X(\partial\Sigma)$ lies in the plane $\{z=z_o\}$, then the surface is unstable for $E_{a,c_o>0,b=0,\alpha,\beta}$.
\end{cor}
{\it Proof.\:} Let $X:\Sigma\longrightarrow{\bf R}^3$ be a critical immersion satisfying \eqref{RME}. Since $b=0$, it follows from \eqref{EL2} and that $z\neq 0$, that $\nu_3=0$ must hold along the boundary $\partial\Sigma$. Then, the variation given by $\delta X=\nu_3\nu$ is a compactly supported normal variation and we can apply Theorem \ref{2E}. 

For $f=\nu_3$ and $b=0$, equation \eqref{2vfE} simplifies to
$$\delta_{\nu_3\nu}^2E[\Sigma]=a\oint_{\partial\Sigma}\frac{\partial_n z}{z}\left(\partial_n \nu_3\right)^2ds\,,$$
since $P[\nu_3]=-2\nu_3/z^2$ holds from \eqref{Pnu3}. Moreover, since $\nu_3$ is a solution of this second order linear differential equation and $\nu_3=0$ on $\partial\Sigma$, we deduce that $\partial_n\nu_3\neq 0$ along $\partial\Sigma$ or, otherwise, $\nu_3=0$ everywhere on $\Sigma$ and the surface would be planar. However, from \eqref{RME} we would obtain that $c_o=-H=0$, contradicting our hypothesis.

Finally, observe that by hypothesis $z>0$ and $\partial_n z<0$ holds on $\partial\Sigma$. Therefore, $\delta_{\nu_3\nu}^2E[\Sigma]<0$ and the surface is unstable. {\bf q.e.d.}

\begin{remark} Of course, the same conclusion follows if $X(\Sigma)$ is contained in $\{z\leq z_o<0\}$ and the boundary $X(\partial\Sigma)$ lies in the plane $\{z=z_o\}$. Indeed, the reflection of a surface across the horizontal plane $\{z=0\}$ preserves the reduced membrane equation \eqref{RME} as well as the quantity $E$.
\end{remark}

The result of Corollary \ref{cor} can be used to show the instability of infinitely many critical domains for $E_{a,c_o>0,b=0,\alpha,\beta}$. We begin by observing that since \eqref{EL1} is satisfied on $\Sigma$, Theorem \ref{rme} and Proposition \ref{converse} also hold in the present case, that is, any smooth axially symmetric topological disc with nonconstant mean curvature critical for $E$ automatically satisfies the reduced membrane equation \eqref{RME}. Furthermore, the topological restriction may be avoided in this case.

\begin{prop}[\cite{PP2}, Proposition 4.2 and Theorem 4.3] Let $X:\Sigma\longrightarrow{\bf R}^3$ be an axially symmetric immersion critical for $E$. Then, either the immersion has constant mean curvature, or else (after a suitable rigid motion and translation of the vertical coordinate) the reduced membrane equation \eqref{RME} holds on $\Sigma$.
\end{prop}

\begin{figure}[h!]
	\makebox[\textwidth][c]{\centering
		\includegraphics[width=3cm]{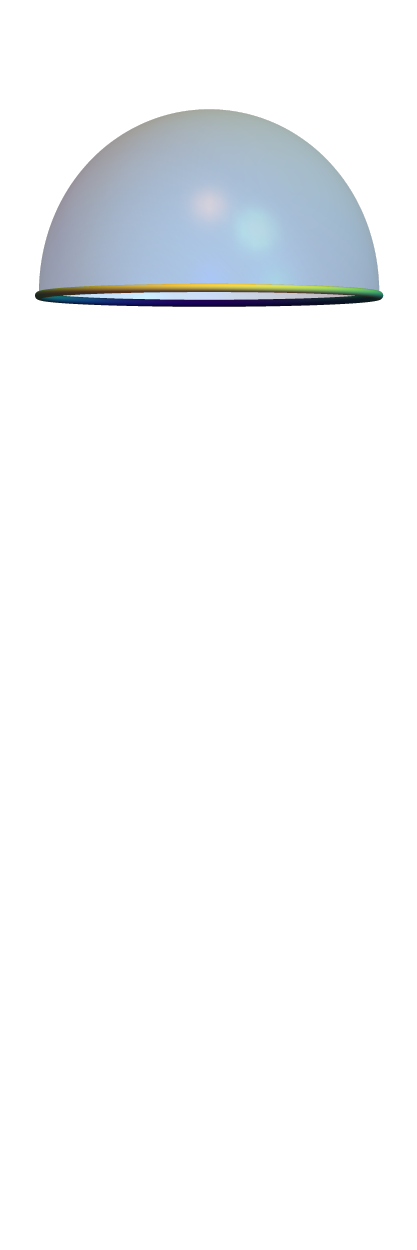}\quad\quad
		\includegraphics[width=3cm]{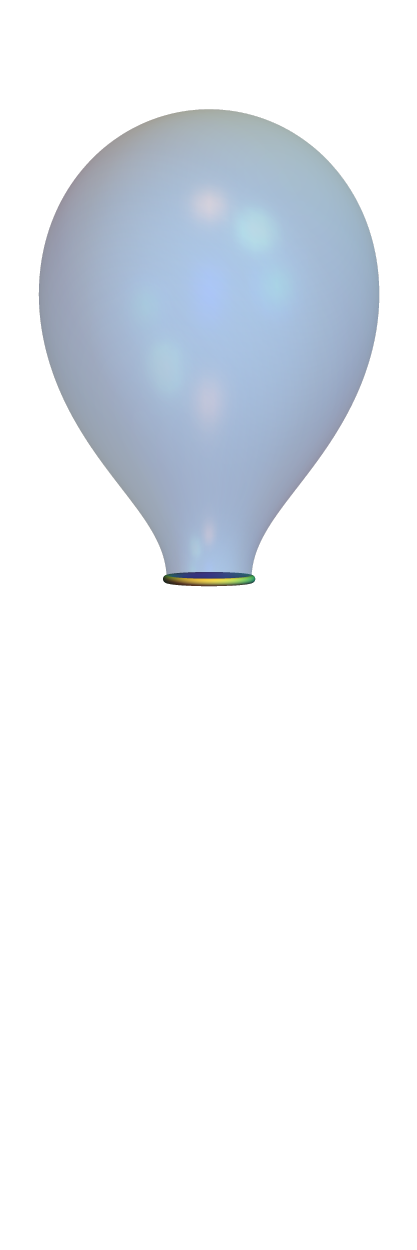}\quad\quad
		\includegraphics[width=3cm]{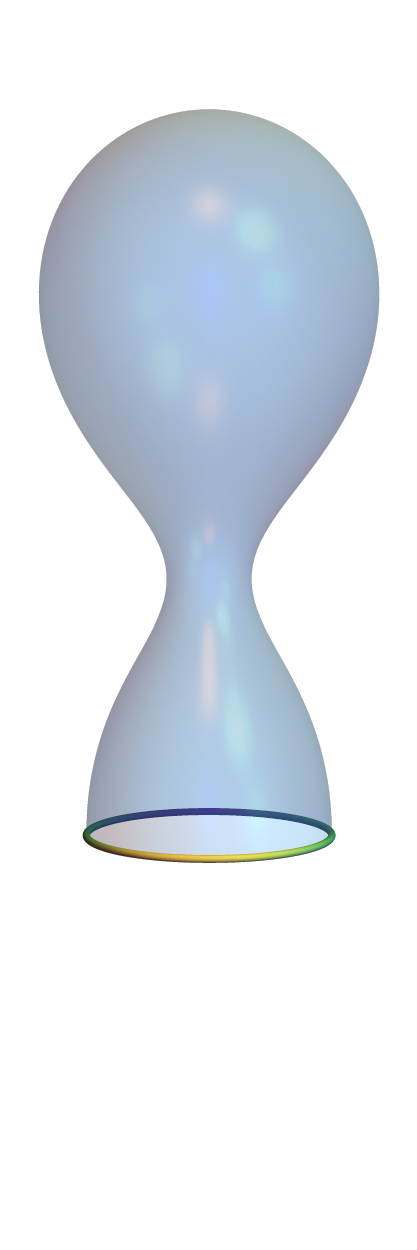}\quad\quad
		\includegraphics[width=3cm]{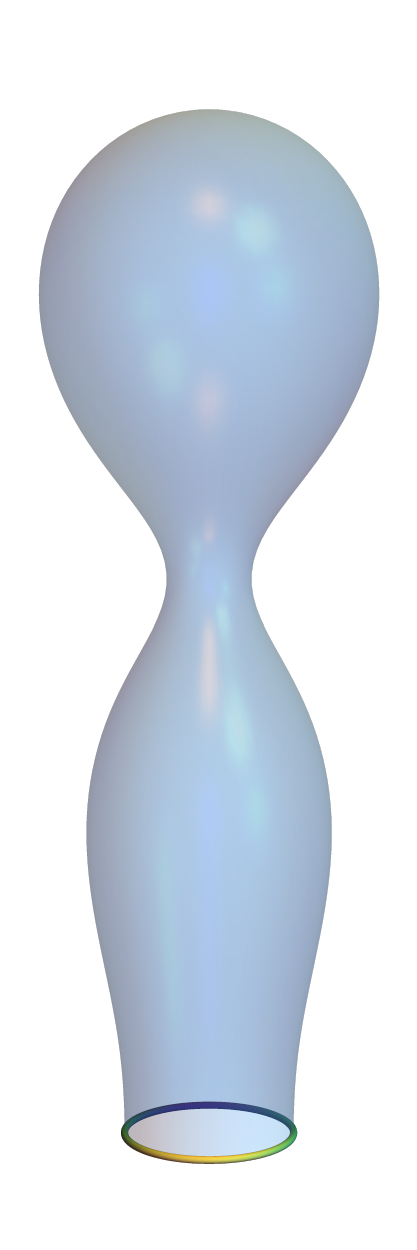}}
	\caption{{\small Four disc type domains in an axially symmetric surface satisfying \eqref{RME}. In all the cases, the boundary is a geodesic circle contained in a horizontal plane located at height $z_o>0$. The axially symmetric surface has been obtained solving the system \eqref{ode1}-\eqref{ode3} with $c_o=2$, initial conditions \eqref{ic} and initial height $\widehat{z}=3$.}}
	\label{EH}
\end{figure}

In order to obtain axially symmetric equilibria for $E$ one can then solve the system of first order ordinary differential equations (c.f. (25)-(27) of \cite{PP2})
\begin{eqnarray}
	r'(\sigma)&=&\cos\varphi(\sigma)\,,\label{ode1}\\
	z'(\sigma)&=&\sin\varphi(\sigma)\,,\label{ode2}\\
	\varphi'(\sigma)&=&-2\frac{\cos\varphi(\sigma)}{z(\sigma)}-\frac{\sin\varphi(\sigma)}{r(\sigma)}-2c_o\,.\label{ode3}
\end{eqnarray}
The solutions of this system represent the arc length parameterized generating curves of axially symmetric surfaces satisfying the reduced membrane equation \eqref{RME}. Moreover, to obtain a regular disc type surface, the initial conditions must be
\begin{equation}\label{ic}
	r(0)=0\,,\quad\quad\quad z(0)=\widehat{z}\neq 0\,,\quad\quad\quad \varphi(0)=0\,,
\end{equation}
where $\widehat{z}\in{\bf R}\setminus\{0\}$ can be assumed to be the parameter. In particular, if this initial height is taken to be positive and sufficiently large it is possible to find a value of $\sigma>0$, say $\sigma_o$, such that $r'(\sigma_o)=0$ and $z\geq z_o=z(\sigma_o)>0$. In other words, it is possible to find an axially symmetric disc type surface $X:\Sigma\longrightarrow{\bf R}^3$ satisfying \eqref{RME} such that the boundary is a geodesic circle and $X(\Sigma)$ is contained in $\{z\geq z_o>0\}$. We point out here that above system \eqref{ode1}-\eqref{ode3} is singular at $r=0$ and, hence, the existence of solution for the initial value problem posed above is not guaranteed by the standard theory of ordinary differential equations. Nevertheless, adapting Proposition 2.2 of \cite{PP3} the existence can be easily proven. Similarly, we highlight that the previous statement about the existence of such value $\sigma_o$ is not clear \emph{a priori}. This will be formally proven in an upcoming work (\cite{LPP}). For the moment and since this is not the objective of the current paper, we just regard this assertion as numerical and we illustrate it in Figure \ref{EH}.

We claim that there are infinitely many pairs of positive constants $(\alpha,\beta)$ such that the domains whose existence has been stated above are critical for an Euler-Helfrich functional $E_{a,c_o>0,b=0,\alpha,\beta}$. Since $b=0$, $\partial\Sigma$ is a geodesic and $z_o>0$, along the boundary $\nu_3=0$ holds and so the boundary conditions \eqref{EL2} and \eqref{EL4} are automatically satisfied. Therefore, the only condition that needs to be checked is \eqref{EL3} which, in this case, reduces to
\begin{equation}\label{Hn}
	a\partial_nH=J'\cdot\nu\,,
\end{equation}
on $\partial\Sigma$.

Denote by $r_o$ the radius of the boundary circle. It follows from the definition of the vector field $J$, see \eqref{Jboundary}, and the derivatives of $T$ along $\partial\Sigma$ (c.f. (2) of \cite{PP1}) that $J'=(\alpha r_o^{-2}-\beta)T'$ and, hence,
$$J'\cdot\nu=\left(\frac{\alpha}{r_o^2}-\beta\right)\kappa_n=-\frac{1}{r_o}\left(\frac{\alpha}{r_o^2}-\beta\right).$$
On the other hand, since \eqref{RME} holds, we have along $\partial\Sigma$,
$$a\partial_nH=a\left(H+c_o\right)_\sigma=-a\left(\frac{\nu_3}{z}\right)_\sigma=a\left(\kappa_1-\left[H+c_o\right]\right)\frac{z_\sigma}{z_o}=-a\kappa_1\frac{1}{z_o}\,,$$
where $\kappa_1$ is the principal curvature such that $\kappa_1=2H-\kappa_n=-2c_o+1/r_o$ along $\partial\Sigma$. Therefore, combining everything and substituting in \eqref{Hn},
$$a\partial_nH=-a\kappa_1\frac{1}{z_o}=a\left(2c_o-\frac{1}{r_o}\right)\frac{1}{z_o}=-\frac{1}{r_o}\left(\frac{\alpha}{r_o^2}-\beta\right).$$
In conclusion, the condition on the pair $(\alpha,\beta)$ is that 
$$a\, \frac{1-2c_o r_o}{z_o}=\alpha r_o^{-2}-\beta\,,$$
holds. Regardless of the value of the left hand side, there are infinitely many positive values of $\alpha$ and $\beta$ solving this equation and thus any of the chosen domains can be considered as a critical surface for an appropriate Euler-Helfrich functional. By Corollary \ref{cor} these surfaces will all be unstable.

The result of Theorem \ref{2E} regarding the second variation of the Euler-Helfrich functional allows us to recover a result of Proposition 6.1 of \cite{PP2}. However, for that result we used a tangential variation of the surface that changed the boundary curve, while in the following corollary we are using a variation arising from vertical translations of the surface.

\begin{cor}\label{cor2}
	Let $X:\Sigma\longrightarrow{\bf R}^3$ be an axially symmetric immersion of a disc type surface critical for $E$ with $b=0$ satisfying the reduced membrane equation and assume that $X(\Sigma)$ is contained in one of the two open half-spaces determined by the plane $\{z=0\}$. If along the boundary circle $\partial\Sigma$
	$$\left(\frac{\alpha}{r_o^2}-\beta\right)K<0\,,$$
	holds, then the surface is unstable for $E_{a,c_o,b=0,\alpha,\beta}$. (Here, $r_o$ denotes the radius of the boundary circle $\partial\Sigma$ and $K$ the Gaussian curvature of the immersion.)
\end{cor} 
{\it Proof.\:} Let $X:\Sigma\longrightarrow{\bf R}^3$ be a critical immersion satisfying \eqref{RME}. Since $b=0$ and $z\neq 0$, the variation $\delta X=\nu_3\nu$ is a compactly supported normal variation. Employing Theorem \ref{2E} for this variation, we have
\begin{eqnarray*}
	\delta_{\nu_3\nu}^2E[\Sigma]&=&a\oint_{\partial\Sigma}\frac{\partial_n z}{z}\left(\partial_n\nu_3\right)^2ds=-a\oint_{\partial\Sigma}\frac{\partial_n\nu_3}{z}\left(2H-\kappa_n\right)ds\\
	&=&a\oint_{\partial\Sigma}\partial_nH\left(2H-\kappa_n\right)ds\,,
\end{eqnarray*}
since $\partial_n\nu_3=-\partial_nz \kappa_1$, where $\kappa_1$ is the principal curvature such that $\kappa_1=2H-\kappa_n$ along $\partial\Sigma$, and
$$\partial_nH=-\partial_n\left(\frac{\nu_3}{z}\right)=-\frac{\partial_n\nu_3}{z}+\nu_3\frac{\partial_n z}{z^2}=-\frac{\partial_n\nu_3}{z}\,,$$
because $\nu_3=0$ holds along $\partial\Sigma$. 

Next, using the Euler-Lagrange equation \eqref{EL3} and the definition \eqref{Jboundary} of $J$, it follows that
$$a\partial_nH=J'\cdot\nu=\left(\frac{\alpha}{r_o^2}-\beta\right)\kappa_n\,.$$
Therefore, combining everything, we obtain
$$\delta_{\nu_3\nu}^2E[\Sigma]=\oint_{\partial\Sigma}\left(\frac{\alpha}{r_o^2}-\beta\right)\kappa_n\left(2H-\kappa_n\right)ds=\oint_{\partial\Sigma}\left(\frac{\alpha}{r_o^2}-\beta\right)K\,ds\,,$$
since, along $\partial\Sigma$, $K=\kappa_n(2H-\kappa_n)$ holds.

Consequently, if the hypothesis of the statement holds, $\delta_{\nu_3\nu}^2E[\Sigma]<0$ and so the surface is unstable. {\bf q.e.d.}

\begin{remark} For both Corollary \ref{cor} and \ref{cor2} we have used the variation $\delta X=\nu_3\nu$. It may seem strange that this variation can destabilize a critical surface since it is the normal component of a variation through vertical translations. However, for the second variation of the Euler-Helfrich functional $E$, the tangential component of the variation field is not insignificant.
\end{remark}

\section{The Operator $P$ as a Jacobi Operator}

From the expression of the second variation of the Helfrich energy given in \eqref{svf2} and Proposition \ref{FP}, it is clear that understanding the second order operator $P$ is essential to analyze the stability or instability of an immersion $X:\Sigma\longrightarrow{\bf R}^3$ satisfying the reduced membrane equation \eqref{RME}. As mentioned in Section 2, the operator $P$ is the linearization of the reduced membrane equation \eqref{RME} for a critical immersion of the functional $\mathcal{G}$. Hence, $P$ can be seen as the Jacobi operator for this functional.

Let $X:\Sigma\longrightarrow{\bf R}^3$ be a smooth immersion of a compact surface $\Sigma$ with smooth boundary and consider the functional $\mathcal{G}$ defined in \eqref{G}. Throughout this section, we will regard the constant $c_o$ as a Lagrange multiplier, so that the critical points of $\mathcal{G}$, viewed as surfaces in the hyperbolic space ${\bf H}^3$, are stationary for the hyperbolic surface area under the constraint that the gravitational potential energy $\mathcal{U}$ is fixed.

Recall that, in terms of the metric in ${\bf R}^3$, the functional $\mathcal{G}$ can be written as
$$\mathcal{G}[\Sigma]=\int_\Sigma \frac{1}{z^2}\,d\Sigma-2c_o\int_V \frac{1}{z^2}\,dV\,.$$
If $X+\epsilon f \nu+\mathcal{O}(\epsilon^2)$, $f\in\mathcal{C}^\infty_o(\Sigma)$, is an admissible normal variation of the immersion $X:\Sigma\longrightarrow{\bf R}^3$, then
\begin{equation}\label{fvfG}
	\delta_{f\nu}\mathcal{G}[\Sigma]=-2\int_\Sigma\left(\xi+c_o\right)\frac{f}{z^2}\,d\Sigma\,,
\end{equation}
where $\xi$ is the function defined in \eqref{xi}. Assume that the immersion $X:\Sigma\longrightarrow{\bf R}^3$ satisfies the reduced membrane equation \eqref{RME}, and so $\xi=-c_o$ holds on $\Sigma$. Then, from the definition of the operator $P$ given in \eqref{P} and \eqref{fvfG}, we easily deduce that the second variation of $\mathcal{G}$ is given by
\begin{equation}\label{svfG}
	\delta_{f\nu}^2\mathcal{G}[\Sigma]=-\int_\Sigma \frac{f P[f]}{z^2}\,d\Sigma\,.
\end{equation}
A critical immersion $X:\Sigma\longrightarrow{\bf R}^3$ is \emph{stable for $\mathcal{G}$} if $\delta_{f\nu}^2\mathcal{G}[\Sigma]\geq 0$ holds for all compactly supported smooth functions $f\in\mathcal{C}_o^\infty(\Sigma)$ satisfying
	\begin{equation}\label{condition}
		\delta_{f\nu}\mathcal{U}[\Sigma]=\int_\Sigma \frac{f}{z^2}\,d\Sigma=0\,,
	\end{equation}
where $\mathcal{U}$ is the gravitational potential energy of $\Sigma$, considered as a surface in ${\bf H}^3$.

\begin{remark}
	A standard argument, using the Implicit Function Theorem, shows that the condition \eqref{condition} is sufficient to produce a compactly supported normal variation of the immersion preserving the value of $\mathcal{U}$.
\end{remark}

We will next show a characterization of the stability condition. In order to state the result, let $\lambda_1<\lambda_2\leq \lambda_3\leq \cdots$ denote the eigenvalues of the Dirichlet boundary problem
\begin{equation}\label{eigenvalueD}
	\left(\frac{P}{z^2}+\lambda\right)[f]=0\,,\quad\quad\quad f\lvert_{\partial\Sigma}=0\,,
\end{equation}
and let $\{f_i\}_{i=1}^\infty$ denote the corresponding eigenfunctions.

\begin{theorem}\label{thmbif}
	Assume that the first two eigenvalues of the problem \eqref{eigenvalueD} satisfy $\lambda_1<0\leq \lambda_2$, and let $X:\Sigma\longrightarrow{\bf R}^3$ be a critical immersion for $\mathcal{G}$. If there exists a solution $h$ of
	\begin{equation}
		P[h]=-2\,,\quad\quad\quad h\lvert_{\partial\Sigma}=0\,,\label{P-2}
	\end{equation}
then the surface is stable if and only if
\begin{equation}\label{stabcond}
	\int_\Sigma \frac{h}{z^2}\,d\Sigma\leq 0\,,
\end{equation}
holds.
\end{theorem}
{\it Proof.\:} Assume that $\lambda_1<0\leq \lambda_2$ holds and that there exists a solution $h$ of the problem \eqref{P-2}.

We will first prove the forward implication. Suppose that the condition \eqref{stabcond} does not hold, that is, suppose
$$\int_\Sigma \frac{h}{z^2}\,d\Sigma>0\,,$$
holds. Since the first eigenfunction $f_1$ of the problem \eqref{eigenvalueD} cannot have a zero on the interior of $\Sigma$, this function may be either positive or negative everywhere on the interior of $\Sigma$. Hence, there exists a constant $\mu\neq 0$, such that
$$\int_\Sigma\frac{\mu f_1+h}{z^2}\,d\Sigma=0\,.$$
Observe that regardless of the sign of the eigenfunction $f_1$, $\mu f_1<0$ holds. 

Moreover, from \eqref{P-2} and the definition of $f_1$ for the problem \eqref{eigenvalueD}, we compute
$$P[\mu f_1+h]=\mu P[f_1]+P[h]=-\mu\lambda_1f_1z^2-2\,.$$
It then follows from this and \eqref{svfG} for $f=\mu f_1+h$, that
\begin{eqnarray*}
	\delta_{f\nu}^2\mathcal{G}[\Sigma]&=&-\int_\Sigma \frac{(\mu f_1+h)P[\mu f_1+h]}{z^2}\,d\Sigma=\int_\Sigma \left(\mu^2\lambda_1f_1^2-\mu h\frac{P[f_1]}{z^2}\right)d\Sigma\\
	&=&\int_\Sigma\left(\mu^2\lambda_1 f_1^2-\mu f_1\frac{P[h]}{z^2}\right)d\Sigma=\int_\Sigma\left(\mu^2\lambda_1f_1^2+2\frac{\mu f_1}{z^2}\right)d\Sigma<0\,.
\end{eqnarray*}
In the third equality above, we have used that the operator $P$ is self-adjoint with respect to the measure $z^{-2}d\Sigma$, while in the last equality we have employed again that $P[h]=-2$ holds. Since $\delta_{f\nu}^2\mathcal{G}[\Sigma]<0$, the surface is unstable.

We next show the converse. Arguing by contradiction, assume that \eqref{stabcond} holds but the surface is unstable. Since the surface is unstable, there exists a compactly supported smooth function $f\in\mathcal{C}_o^\infty(\Sigma)$ such that \eqref{condition} holds and
\begin{equation}\label{help}
	\delta_{f\nu}^2\mathcal{G}[\Sigma]=-\int_\Sigma \frac{f P[f]}{z^2}\,d\Sigma<0\,.
\end{equation} 
We can find constants $\mu,\eta\in{\bf R}$ (not both zero), such that
$$\int_\Sigma f_1\left(\mu h+\eta f\right)d\Sigma=0\,,$$
where $f_1$ is the first eigenfunction of the problem \eqref{eigenvalueD}.

Consider first the possibility that $\eta=0$. Since $\eta=0$, $\mu$ must be different from zero. Then, from the eigenvalue problem \eqref{eigenvalueD}, we would obtain
$$0=\int_\Sigma f_1h\,d\Sigma=\frac{-1}{\lambda_1}\int_\Sigma h\frac{P[f_1]}{z^2}\,d\Sigma=\frac{-1}{\lambda_1}\int_\Sigma f_1\frac{P[h]}{z^2}\,d\Sigma=\frac{2}{\lambda_1}\int_\Sigma\frac{f_1}{z^2}\,d\Sigma\,,$$
which contradicts the fact that $f_1$ cannot have zeros on the interior of $\Sigma$. Therefore, $\eta\neq 0$ holds necessarily.

We then compute for $g=\mu h+\eta f$,
\begin{eqnarray*}
	\delta_{g\nu}^2\mathcal{G}[\Sigma]&=&\int_\Sigma-\frac{(\mu h+\eta f)P[\mu h+\eta f]}{z^2}\,d\Sigma\\
	&=&\int_\Sigma\left(-\mu^2 \frac{h P[h]}{z^2}-2\mu \eta \frac{f P[h]}{z^2}-\eta^2 \frac{f P[f]}{z^2}\right)d\Sigma\\
	&=&\int_\Sigma\left(2\mu^2\frac{h}{z^2}+4\mu\eta\frac{f}{z^2}-\eta^2\frac{f P[f]}{z^2}\right)d\Sigma\\
	&=&\int_\Sigma \left(2\mu^2\frac{h}{z^2}-\eta^2\frac{fP[f]}{z^2}\right)d\Sigma<0\,,
\end{eqnarray*}
where in the third line we have used that $P[h]=-2$ holds and in the fourth one that \eqref{condition} holds for the function $f$. The above inequality, which follows from \eqref{stabcond} and \eqref{help}, contradicts the hypothesis that $\lambda_2\geq 0$ holds, proving the result. {\bf q.e.d.}

\begin{remark} If $\lambda_2>0$ holds, the function $h$ solution of \eqref{P-2} always exists. If, instead, $\lambda_2=0$ holds, then there is still a characterization of stability analogous to that for the constant mean curvature case (see \cite{K}). For example, it may occur that no solution of \eqref{P-2} exists and, in this case, the surface is unstable.
\end{remark}

\begin{figure}[h!]
	\makebox[\textwidth][c]{\centering
		\includegraphics[height=4cm]{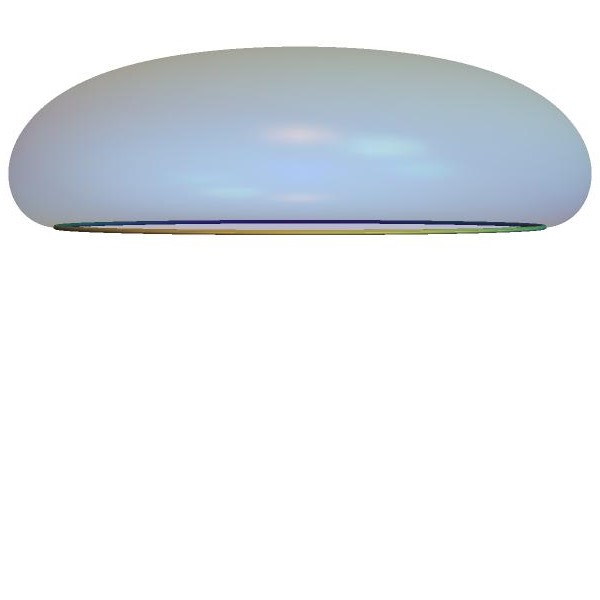}\quad\quad
		\includegraphics[height=4cm]{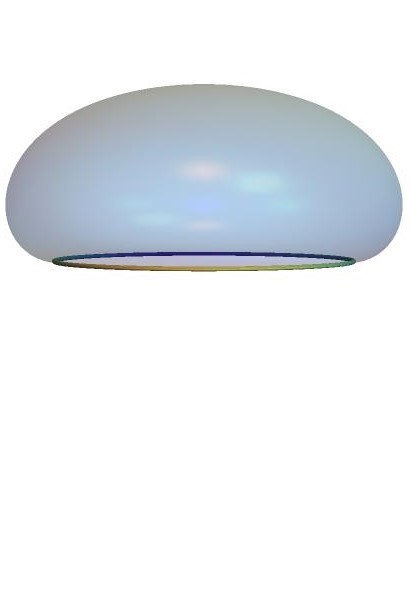}\quad\quad
		\includegraphics[height=4cm]{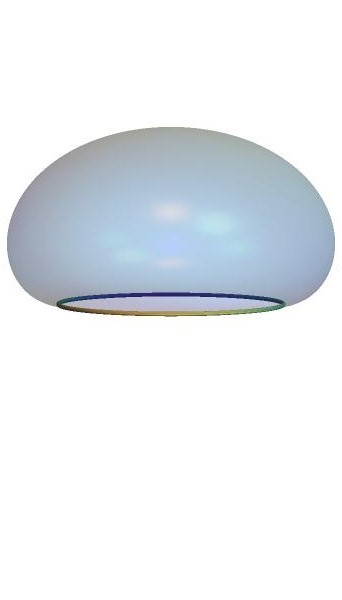}\quad\quad
		\includegraphics[height=4cm]{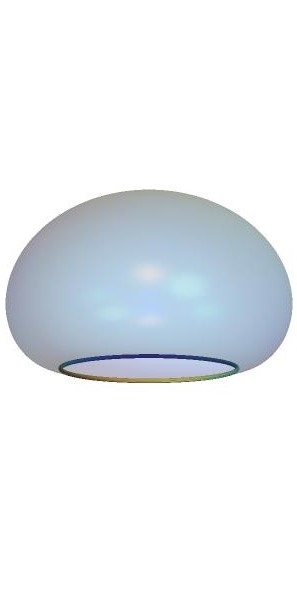}}
	\caption{{\small Four axially symmetric convex embedded topological discs satisfying \eqref{RME} and such that, along their boundary circles, the tangent plane is horizontal. They have been obtained solving the system \eqref{ode1}-\eqref{ode3} with $c_o=2$, initial conditions \eqref{ic} and initial heights $\widehat{z}=-0.55$, $\widehat{z}=-0.7$, $\widehat{z}=-0.9$, and $\widehat{z}=-1.2$, respectively.}}
	\label{Bifurcation}
\end{figure}

In what follows we will apply the result of Theorem \ref{thmbif} to a class of critical immersions for $\mathcal{G}$ discussed in \cite{PP3}. The immersions of these surfaces, denoted by $\Sigma_0$, are axially symmetric convex embedded topological discs satisfying the reduced membrane equation \eqref{RME} and having a horizontal tangent plane along its boundary circle. In addition, possibly after a reflection across the plane $\{z=0\}$, $X(\Sigma_0)$ are contained in $\{z\leq \widehat{z}<0\}$, where $\widehat{z}$ is their maximum height. Some of these surfaces are illustrated in Figure \ref{Bifurcation}. The geometric description of these surfaces, which arise for $\widehat{z}<-1/c_o<0$, was formally proven in Theorem 2.2 of \cite{PP3}.

We will show that the immersions of these surfaces $\Sigma_0$ are stable for $\mathcal{G}$.

\begin{cor}\label{stable} Axially symmetric convex embedded topological discs satisfying the reduced membrane equation \eqref{RME} and having horizontal planes along their boundaries are stable for $\mathcal{G}$.
\end{cor}
{\it Proof.\:} Let $\Sigma_0$ be a topological disc and $X:\Sigma_0\longrightarrow{\bf R}^3$ an axially symmetric convex embedding of such a surface satisfying the reduced membrane equation \eqref{RME} and having horizontal tangent plane along $X(\partial\Sigma_0)$. In this case, the solution of \eqref{P-2} is (\cite{PP3})
$$h:=c_o^{-1}\left(q(\ell)\psi-\psi(\ell)q\right),$$
where $\psi$ is a non trivial axially symmetric solution of $P[\psi]=0$ (its existence is shown in Lemma 3.1 of \cite{PP3}) and $q:=X\cdot \nu$ is the support function of $\Sigma_0$. Here, $\sigma=\ell$ represents the boundary of the surface (notice that in \cite{PP3} the boundary is taken to be $\sigma=0$). For simplicity, we will take that the solution $\psi$ is normalized by $\psi(\ell)=1$ (recall that from Lemma 3.2 of \cite{PP3}, $\psi(\ell)\neq 0$). 

Since the problem \eqref{P-2} has a solution, from Theorem \ref{thmbif}, it suffices to show that the condition \eqref{stabcond} holds in order to prove that the surface $X(\Sigma_0)$ is stable. From the definition of the solution $h$ given above, we have that
\begin{equation}\label{split}
	\int_{\Sigma_0}\frac{h}{z^2}\,d\Sigma=c_o^{-1}\left(q(\ell)\int_{\Sigma_0}\frac{\psi}{z^2}\,d\Sigma-\int_{\Sigma_0}\frac{q}{z^2}\,d\Sigma\right),
\end{equation}
since we are taking $\psi(\ell)=1$. Observe that since $X(\Sigma_0)$ is contained in $\{z\leq \widehat{z}<0\}$ the normal $\nu$ to the surface along the boundary is the vector $(0,0,-1)$. Hence, $q(\ell)=-z(\ell)>-\widehat{z}>0$ holds. Moreover, we also notice here that in \cite{PP3}, using an argument involving the rescaling of the immersion $X$, it was shown that $P[q]=2c_o$ holds. Combining both things together with $P[\psi]=0$, we obtain an upper bound of the first term on the right hand side of \eqref{split}. In fact,
\begin{eqnarray*}
	\int_{\Sigma_0}\frac{\psi}{z^2}\,d\Sigma&=&\frac{1}{2c_o}\int_{\Sigma_0}\psi\frac{P[q]}{z^2}\,d\Sigma=\frac{1}{2c_o}\int_{\Sigma_0}\left(\psi\frac{P[q]}{z^2}-q\frac{P[\psi]}{z^2}\right)d\Sigma\\
	&=&\frac{1}{2c_o}\int_{\Sigma_0}\left(\psi \nabla\cdot\left[\frac{\nabla q}{z^2}\right]-q \nabla\cdot\left[\frac{\nabla \psi}{z^2}\right]\right)d\Sigma\\
	&=&\frac{1}{2c_o}\oint_{\partial\Sigma_0}\left(\frac{\partial_nq}{z^2}-q(\ell)\frac{\partial_n\psi}{z^2}\right)ds\,,
\end{eqnarray*}
where in the last equality we have used Green's second identity and that $\psi(\ell)=1$. Since $q(\ell)>0$ and $\partial_n\psi(\ell)>0$ holds along $\partial\Sigma_0$ (for the latter assertion, see the proof of Lemma 4.1 of \cite{PP3}), it follows that the second term in above boundary integral is negative and, hence,
$$\int_{\Sigma_0}\frac{\psi}{z^2}\,d\Sigma<\frac{1}{2c_o}\oint_{\partial\Sigma_0}\frac{\partial_nq}{z^2}\,ds=\frac{\pi\, r(\ell)\partial_nq(\ell)}{c_o z^2(\ell)}\,.$$
We now compute
\begin{eqnarray*}
	\partial_nq(\ell)&=&\partial_n\left(X\cdot \nu\right)(\ell)=X\cdot d\nu(n)(\ell)=\kappa(\ell)X\cdot n(\ell)=\kappa(\ell)r(\ell)r'(\ell)\\&=&-\varphi'(\ell)r(\ell)r'(\ell)
	=2\frac{r(\ell)}{z(\ell)}-2c_or(\ell)\,,
\end{eqnarray*}
where in the last equality we have used the system \eqref{ode1}-\eqref{ode2} and that $\varphi(\ell)=-\pi$. Therefore, since $c_o^{-1}q(\ell)>0$, it follows that the first term in the right hand side of \eqref{split} satisfies
\begin{equation}\label{help1}
	c_o^{-1}q(\ell)\int_{\Sigma_0}\frac{\psi}{z^2}\,d\Sigma<-2\pi\, \frac{r^2(\ell)}{c_o^2 z^2(\ell)}+2\pi\,\frac{r^2(\ell)}{c_o z(\ell)}\,.
\end{equation}

We next work with the second term. For this integral, we apply the Divergence Theorem to
$$\int_V \frac{1}{z^2}\,dV=\int_V\nabla\cdot\left(\frac{X}{z^2}\right)dV=\int_{\Sigma_0}\frac{q}{z^2}\,d\Sigma+\int_D\frac{-1}{z}\,d\Sigma\,,$$
where $V$ is the region bounded by $\Sigma_0$ and the planar disc $D$ at its base (which is located at height $z(\ell)$ and has radius $r(\ell)$). Consequently,
\begin{equation}\label{help2}
	\int_{\Sigma_0}\frac{q}{z^2}\,d\Sigma=\int_V\frac{1}{z^2}\,dV+\int_D\frac{1}{z}\,d\Sigma>\int_D\frac{1}{z}\,d\Sigma=\pi\,\frac{r^2(\ell)}{z(\ell)}\,.
\end{equation}

Finally, using \eqref{help1} and \eqref{help2} in \eqref{split}, we deduce that
\begin{equation*}
	\int_{\Sigma_0}\frac{h}{z^2}\,d\Sigma<-2\pi\frac{r^2(\ell)}{c_o^2 z^2(\ell)}+2\pi\frac{r^2(\ell)}{c_o z(\ell)}-\pi\frac{r^2(\ell)}{c_oz(\ell)}=-2\pi\frac{r^2(\ell)}{c_o^2 z^2(\ell)}+\pi\frac{r^2(\ell)}{c_o z(\ell)},
\end{equation*}
which is clearly negative since $z(\ell)<0$ holds. This finishes the proof. {\bf q.e.d.}
\\

As a consequence of Corollary \ref{stable} and by monotonicity, we have that any domain $\Omega\subset\Sigma_0$ is also stable for $\mathcal{G}$. We refer to these domains $\Omega$ as \emph{subdomains of $\Sigma_0$}. An illustration of such a subdomain is the surface (A) of Figure 2 of \cite{PP3}. By the same reasoning, any domain $\widehat{\Omega}$ strictly containing $\Sigma_0$ is unstable for $\mathcal{G}$. The domains $\widehat{\Omega}$ are called \emph{superdomains of $\Sigma_0$} and a couple of examples are the surfaces (C) and (D) of Figure 2 of \cite{PP3}. We summarize this in the following result.

\begin{cor}
	Subdomains of $\Sigma_0$ are stable and superdomains of $\Sigma_0$ are unstable for the functional 
	$\mathcal{G}$.
\end{cor}

\begin{remark} In Theorem 3.1 of \cite{PP3} it was shown that the surface $\Sigma_0$ is embedded in a one parameter family of axially symmetric solutions of the reduced membrane equation \eqref{RME}, parameterized by $c_o$, which all share the same boundary circle. At the critical value of the parameter at which the tangent planes along the boundary circle are horizontal (i.e., precisely at $\Sigma_0$), a non-axially symmetric branch of solutions of \eqref{RME} bifurcates while maintaining the same circular boundary. The new branch of surfaces has a vertical plane of reflective symmetry and this branch together with its reflected image constitute a pitchfork bifurcation, that is in all likelihood subcritical.

Again, by virtue of the fact that all these surfaces satisfy the reduced membrane equation \eqref{RME}, they automatically satisfy \eqref{EL} as well. 
\end{remark}

\bigskip

\begin{flushleft}
	Bennett P{\footnotesize ALMER}\\
	Department of Mathematics,
	Idaho State University,
	Pocatello, ID 83209,
	U.S.A.\\
	E-mail: palmbenn@isu.edu
\end{flushleft}

\bigskip

\begin{flushleft}
	\'Alvaro P{\footnotesize \'AMPANO}\\
	Department of Mathematics and Statistics, Texas Tech University, Lubbock, TX
	79409, U.S.A.\\
	E-mail: alvaro.pampano@ttu.edu
\end{flushleft}

\end{document}